\documentclass[11pt]{article}
\usepackage{latexsym,amssymb}
\def\demo{\noindent{\bf Proof. }}
\def\QED{\hfill$\Box$}
\newtheorem{Theorem}{Theorem}[section]
\newtheorem{Lemma}[Theorem]{Lemma}
\newtheorem{Corollary}[Theorem]{Corollary}
\newtheorem{Proposition}[Theorem]{Proposition}
\newtheorem{Remark}[Theorem]{Remark}
\newtheorem{Example}[Theorem]{Example}
\newtheorem{Conjecture}[Theorem]{Conjecture}
\newtheorem{Definition}[Theorem]{Definition}

\begin{document}
\topmargin3mm
\hoffset=-1cm
\voffset=-1.5cm
\

\medskip

\begin{center}
{\large\bf Blowup algebras of square-free monomial ideals and some
links   
\\  to combinatorial optimization problems}
\vspace{6mm}\\

\medskip

Isidoro Gitler,\ \,  
Enrique Reyes\footnote{Partially supported by COFAA-IPN.}, 
 \\ 
and 
\\ 
Rafael H. Villarreal\footnote{Partially supported by CONACyT 
grant 49251-F and SNI.} 
\\ 
{\small Departamento de Matem\'aticas}\vspace{-1mm}\\ 
{\small Centro de Investigaci\'on y de Estudios Avanzados del
IPN}\vspace{-1mm}\\   
{\small Apartado Postal 14--740}\vspace{-1mm}\\ 
{\small 07000 M\'exico City, D.F.}\vspace{-1mm}\\ 
{\small e-mail: {\tt vila@math.cinvestav.mx}}\vspace{4mm}
\end{center}
\date{}

\begin{abstract} 
\noindent Let $I=(x^{v_1},\ldots,x^{v_q})$ be a 
square-free monomial ideal of a polynomial ring $K[x_1,\ldots,x_n]$ 
over an arbitrary field $K$ and let $A$ be the incidence matrix with
column vectors  
${v_1},\ldots,{v_q}$.  We will establish some connections between
algebraic  
properties of certain graded algebras associated to $I$ and
combinatorial  
optimization properties of certain polyhedra and clutters  
associated to $A$ and $I$ respectively. Some applications to Rees 
algebras and combinatorial optimization are presented. We study a
conjecture of Conforti and 
Cornu\'ejols using an algebraic approach.
\end{abstract}

\section{Introduction}\label{Int}

Let $R=K[x_1,\ldots,x_n]$ be a polynomial ring 
over a field $K$ and let $I$ be an ideal 
of $R$ of height $g\geq 2$ minimally generated by a finite set of
square-free monomials  
$F=\{x^{v_1},\ldots,x^{v_q}\}$ of degree at least two. As usual we
use  
$x^a$ as an abbreviation for $x_1^{a_1} \cdots x_n^{a_n}$, 
where $a=(a_1,\ldots,a_n)\in \mathbb{N}^n$. A {\it clutter\/} with
vertex set  
$X$ is a family of subsets of $X$, called edges, none 
of which is included in another. We associate to the 
ideal $I$ a {\it clutter\/} $\cal C$ by taking the set 
of indeterminates $X=\{x_1,\ldots,x_n\}$ as vertex set and 
$E=\{S_1,\ldots,S_q\}$ as edge set, where 
$$S_k={\rm supp}(x^{v_k})=\{x_i\vert\, \langle e_i,v_k\rangle=1\}.$$
Here $\langle\ ,\, \rangle$ denotes the standard 
inner product and $e_i$ is the $i${\it th} unit vector. For this
reason  
$I$ is called the {\it edge ideal\/} of $\cal C$. To stress the
relationship  
between $I$ and $\cal C$ we will use the notation $I=I({\cal C})$. A
basic example  
of clutter is a graph. Algebraic 
and combinatorial properties 
of edge ideals and graded algebras associated to graphs have been  
studied in \cite{susan1,herzog-survey,ITG,carlos-tesis,Vi2}. The
related 
notion of facet ideal has been studied by Faridi 
\cite{faridi,faridi1} and Zheng \cite{zheng}.  

The {\it blowup algebras\/} studied here are 
the {\it Rees algebra\/}  
$$
R[It]=R\oplus It\oplus\cdots\oplus I^{i}t^i\oplus\cdots
\subset R[t],
$$
where $t$ is a new variable, and the {\it associated graded 
ring\/}
$$
{\rm gr}_I(R)=R/I\oplus I/I^2\oplus\cdots\oplus 
I^i/I^{i+1}\oplus\cdots\simeq R[It]\otimes_R(R/I),
$$
with multiplication  
$(a+I^{i+1})(b+I^{j+1})=ab+I^{i+j+1}$, $a\in I^{i}$, $b\in I^{j}$. 

In the sequel $A$ will denote the {\it incidence matrix\/} of order
$n\times q$ whose column vectors are 
$v_1,\ldots,v_q$. In order to link the properties of these algebras
with combinatorial  
optimization problems we consider 
the {\it set covering polyhedron\/}
$$
Q(A)=\{x\in\mathbb{R}^n\vert\, x\geq 0;\, xA\geq\mathbf{1}\}, 
$$
and the related system of linear inequalities $x\geq 0;\ xA\geq\mathbf{1}$, 
where $\mathbf{1}=(1,\ldots,1)$. Recall that this system is called 
{\it totally dual integral\/} (TDI) if the maximum in 
the LP-duality equation
\begin{equation}\label{jun6-2-03}
{\rm min}\{\langle \alpha,x\rangle \vert\, x\geq 0; xA\geq \mathbf{1}\}=
{\rm max}\{\langle y,\mathbf{1}\rangle \vert\, y\geq 0; Ay\leq\alpha\} 
\end{equation}
has an integral optimum solution $y$ for each integral vector $\alpha$ with 
finite maximum. If the system is totally dual integral 
it is seen that $Q(A)$ 
has only integral vertices, this follows from 
\cite[Theorem~22.1, Corollary~22.1.a, pp.~310-311]{Schr}.

We are able to express algebraic properties of blowup algebras in terms 
of TDI systems 
and combinatorial properties of clutters, such as the integrality 
of $Q(A)$ and the K\"onig property. An important goal here is to 
establish bridges between commutative algebra and combinatorial 
optimization, which 
could be beneficial to both areas. Necessary and/or 
sufficient conditions for the normality of $R[It]$ and the 
reducedness of ${\rm gr}_I(R)$ are shown. Some of our results give
some support  
to a conjecture of Conforti and 
Cornu\'ejols (Conjecture~\ref{conforti-cornuejols1}). Applications to
Rees algebras  
theory and combinatorial optimization are presented. 

Along the paper we introduce some of the algebraic and combinatorial 
notions that are most relevant. For unexplained terminology and
notation we refer to  
\cite{korte,oxley,Schr} and \cite{Mats,Vas}. See \cite{cornu-book}
for detailed  
information about clutters.  

\section{Vertex covers of clutters} 

The set of non-negative real numbers will 
be denoted by $\mathbb{R}_+$. 
To avoid repetitions throughout this article we shall use the notation 
and assumptions introduced in Section~\ref{Int}. For convenience we
shall always assume that each variable $x_i$ occurs in at least one
monomial of $F$. 

\begin{Definition}\rm A subset $C\subset X$ is a 
{\it minimal vertex cover\/} of the clutter $\cal C$ if: 
(i) every edge of $\cal C$ contains at least one vertex of $C$, 
and (ii) there is no proper subset of $C$ with the first 
property. If $C$ satisfies condition (i) only, then $C$ is 
called a {\it vertex cover\/} of $\cal C$.   
\end{Definition}

The first aim is to characterize this notion in terms 
of the integral vertices of set covering polyhedrons and the minimal 
primes of edge ideals. 

\medskip

\noindent {\it Notation\/} The {\it support\/} of 
$x^a=x_1^{a_1}\cdots x_n^{a_n}$ is ${\rm supp}(x^a)= \{x_i\, |\,
a_i>0\}$.

\begin{Proposition}\label{1cover} The following are equivalent\/{\rm :}
\begin{description}
\item{\rm(a)} $\mathfrak{p}=(x_{1},\ldots,x_{r})$ is a minimal prime of 
$I=I({\cal C})$. 
\item{\rm(b)} $C=\{x_{1},\ldots,x_{r}\}$ is a minimal vertex cover of
$\cal C$. 
\item{\rm(c)} $\alpha=e_{1}+\cdots+e_{r}$ is a vertex of $Q(A)$.
\end{description}
\end{Proposition}

\demo (a) $\Leftrightarrow$ (b): It follows readily by noticing that
the minimal  
primes of the square-free 
monomial ideal $I$ are face 
ideals, that is, they are generated by 
subsets of the set of variables, see \cite[Proposition~5.1.3]{monalg}.  

(b) $\Rightarrow$ (c): Fix $1\leq i\leq r$. To make notation simpler
fix $i=1$. We  
may assume that there is an $s_1$ such that $x^{v_j}=x_1m_j$ for
$j=1,\ldots,s_1$  
and $x_1\notin{\rm supp}(x^{v_j})$ for $j>s_1$. Notice that 
${\rm supp}(m_{k_1})\cap(C\setminus\{x_1\})=\emptyset$ for 
some $1\leq k_1\leq s_1$, otherwise $C\setminus\{x_1\}$ is a vertex
cover of $\cal C$  
strictly contained in $C$, a contradiction. Thus 
${\rm supp}(m_{k_1})\cap C=\emptyset$ because $I$ is square-free.
Hence for each $1\leq 
i\leq r$ there is $v_{k_i}$ in $\{v_1,\ldots,v_q\}$ such that
$x^{v_{k_i}}=x_im_{k_i}$ and  
${\rm supp}(m_{k_i})\subset\{x_{r+1},\ldots,x_n\}$. The vector $\alpha$ 
is clearly in $Q(A)$, and since 
$\{e_i\}_{i=r+1}^n\cup\{v_{k_1},\ldots,v_{k_r}\}$ is linearly
independent, and  
$$
\langle\alpha,e_i\rangle=0\ \ \ \ (i=r+1,\ldots,n); 
\ \ \ \langle\alpha,v_{k_i}\rangle=1\ \ \ \ (i=1,\ldots,r),
$$
we get that the vector $\alpha$
is a basic feasible solution. Therefore by
\cite[Theorem~2.3]{bertsimas} 
$\alpha$ is a 
vertex of $Q(A)$. 

(c) $\Rightarrow$ (b): It is clear that $C$ intersects all the edges
of the clutter $\cal C$ because $\alpha\in Q(A)$. If $C'\subsetneq C$
is a vertex cover of $\cal C$, then  
the vector $\alpha'=\sum_{x_i\in C'}e_i$ satisfies $\alpha'A\geq
\mathbf{1}$ and $\alpha'\geq 0$.  
Using that $\alpha$ is a basic feasible solution in the sense of
\cite{bertsimas} it is not  
hard to verify that $\alpha'$ is also a vertex of $Q(A)$. By the finite 
basis theorem \cite[Theorem 4.1.3]{webster} we can write
$$
Q(A)=\mathbb{R}_+^n+{\rm conv}(V),
$$
where $V$ is the vertex set of $Q(A)$. As $\alpha=\beta+\alpha'$, for some 
$0\neq \beta\in\mathbb{R}_+^n$, we get 
$$
Q(A)=\mathbb{R}_+^n+{\rm conv}(V\setminus\{\alpha\}).
$$
Hence the vertices of $Q(A)$ are contained in $V\setminus\{\alpha\}$
(see \cite[Theorem~7.2]{Bron}), 
a contradiction. 
Thus $C$ is a 
minimal vertex cover. \QED

\begin{Corollary}\label{may12-06} A vector $\alpha\in\mathbb{R}^n$ 
is an integral vertex of $Q(A)$  if and only if $\alpha$ is equal to
$e_{i_1}+\cdots+e_{i_s}$ for  
some minimal vertex cover $\{x_{i_1},\ldots,x_{i_s}\}$ of $\cal C$. 
\end{Corollary}

\demo By Proposition~\ref{1cover} it suffices to observe that any
integral vertex of $Q(A)$ has entries in $\{0,1\}$ because $A$ has
entries in $\{0,1\}$. See \cite[Lemma~4.6]{shiftcon}. \QED

\medskip

A set of edges of the clutter $\cal C$ is {\it independent} if no two
of them have a common vertex.  
We denote the smallest number of vertices in any 
minimal vertex cover of $\cal C$ by ${\alpha}_0({\cal C})$ and the 
maximum number of independent edges of ${\cal C}$ by $\beta_1({\cal
C})$. These  
numbers are related to min-max problems because they satisfy:
\begin{eqnarray*}
\lefteqn{\alpha_0({\cal C})\geq {\rm min}\{\langle
\mathbf{1},x\rangle \vert\, x\geq 0; xA\geq  
\mathbf{1}\}}\\
&\ \ \ \ \ \ \ \ \ \ &
={\rm max}\{\langle y,\mathbf{1}\rangle \vert\, y\geq 0;
Ay\leq\mathbf{1}\} 
\geq \beta_1({\cal C}). 
\end{eqnarray*}
Notice that $\alpha_0({\cal C})=\beta_1({\cal C})$ if and only if
both sides of  
the equality have integral optimum solutions. 

These two numbers can be interpreted in terms of $I$. By
Proposition~\ref{1cover} the height  
of the ideal $I$, denoted by ${\rm ht}(I)$, is equal to the {\it covering
number\/} $\alpha_0({\cal C})$. On the 
other hand the {\it independence number\/} $\beta_1({\cal C})$ is
equal to ${\rm mgrade}(I)$,  
the {\it monomial grade} of the ideal:
$$
\beta_1({\cal C})=\max\{r\vert\, \exists\, \mbox{ a regular sequence
of monomials } 
x^{\alpha_1},\ldots,x^{\alpha_r}\in I\}.
$$
The equality $\alpha_0({\cal C})=\beta_1({\cal C})$ is equivalent 
to require $x_1\cdots x_nt^g\in R[It]$, where $g$ is the 
covering number $\alpha_0({\cal C})$.

\begin{Definition}\rm If $\alpha_0({\cal C})=\beta_1({\cal C})$ we
say that the clutter  
$\cal C$ (or the ideal $I$) has the {\it K\"onig property\/}.   
\end{Definition}

\section{Rees algebras and polyhedral geometry}
\label{valencia}

Let ${\cal A}=\{v_1,\ldots,v_q\}$ be the set of exponent vectors of 
$x^{v_1},\ldots,x^{v_q}$ and let  
$$
{\cal
A}'=\{e_1,\ldots,e_n,(v_1,1),\ldots,(v_q,1)\}\subset\mathbb{R}^{n+1}, 
$$
where $e_i$ is the $i${\it th} unit vector. The {\it Rees cone\/} of 
${\cal A}$ is the rational 
polyhedral cone, denoted by ${\mathbb R}_+{\cal A}'$, consisting of the 
linear combinations of ${\cal A}'$ with non-negative 
coefficients. Note $\dim(\mathbb{R}_+{\cal A}')=n+1$. Thus according 
to \cite{webster} there is a unique irreducible representation 
\begin{equation}\label{okayama} 
{\mathbb R}_+{\cal A}'=H_{e_1}^+\cap \cdots\cap H_{e_{n+1}}^+
\cap H_{a_1}^+\cap\cdots\cap H_{a_r}^+
\end{equation}
such that $0\neq a_i\in\mathbb{Q}^{n+1}$ and 
$\langle a_i,e_{n+1}\rangle=-1$ for all 
$i$. As usual $H_{a}^+$ denotes 
the closed halfspace
$$
H_a^+=\{\alpha\in\mathbb{R}^{n+1}\vert\, \langle
\alpha,a\rangle\geq 0\}
$$
and $H_a$ is the hyperplane through the origin with normal vector $a$. 

\begin{Theorem}\label{gvvalencia-theo} The function 
$\varphi\colon \mathbb{Q}^n\rightarrow\mathbb{Q}^{n+1}$ 
given by $\varphi(\alpha)=(\alpha,-1)$ induces a bijective map 
$$
\varphi\colon V\longrightarrow\{a_1,\ldots,a_r\}
$$
between the set of vertices $V$ of $Q(A)$ and the set
$\{a_1,\ldots,a_r\}$  
of normal vectors  that occur in the irreducible representation of\/
$\mathbb{R}_+{\cal A}'$. 
\end{Theorem}

\demo First we show the containment
$\varphi(V)\subset\{a_1,\ldots,a_r\}$.  
Take $\alpha$ in $V$\/. By \cite[Theorem~2.3]{bertsimas} $\alpha$  
is a basic feasible solution. Hence 
$\langle \alpha,v_i\rangle\geq 1$ for $i=1,\ldots,q$, $\alpha\geq 0$,
and  
there exist $n$ linearly independent vectors $v_{i_1},\ldots,v_{i_k}, 
e_{j_1},\ldots,e_{j_s}$  in ${\cal A}\cup\{e_1,\ldots,e_n\}$ such
that  
$\langle \alpha,v_{i_{h}}\rangle=1$ and $\langle
\alpha,e_{j_m}\rangle=0$  
for all ${h},m$. It follows that the set 
$$F=H_{(\alpha,-1)}\cap\mathbb{R}_+{\cal A}'$$ 
has  dimension $n$ and $\mathbb{R}_+{\cal A}'\subset
H_{(\alpha,-1)}^+$.  
Therefore $F$ is a facet of $\mathbb{R}_+{\cal A}'$. Using
\cite[Theorem~3.2.1]{webster}  
we obtain that $F=\mathbb{R}_+{\cal A}'\cap H_{a_p}$ for some 
$1\leq p\leq r$, and  
consequently $H_{(\alpha,-1)}=H_{a_p}$. Since the first $n$ entries
of $a_p$ are  
non-negative and $\langle a_p,e_{n+1}\rangle=-1$ it follows 
that $\varphi(\alpha)=(\alpha,-1)=a_p$, as desired.  

To show the reverse containment write $a_p=(\alpha,-1)$, with $1\leq
p\leq r$ and  
$\alpha\in\mathbb{R}^n$. We will prove that $\alpha$ is a vertex of
$Q(A)$.  
Since Eq.~(\ref{okayama}) is an irreducible representation one has
that the set 
$$
F=H_{(\alpha,-1)}\cap\mathbb{R}_+{\cal A}'
$$ 
is a facet of the Rees cone $\mathbb{R}_+{\cal A}'$, 
see \cite[Theorem~3.2.1]{webster}. Hence there is a 
linearly independent set 
$$
\{(v_{i_1},1),\ldots,(v_{i_k},1),e_{j_1},\ldots,e_{j_s}\}\subset{\cal
A}'\ \ \ \ (k+s=n) 
$$
such that
\begin{eqnarray}
\langle(v_{i_{h}},1),(\alpha,-1)\rangle=0&\Rightarrow&
\langle v_{i_{h}},\alpha\rangle=1\ \ \ \
({h}=1,\ldots,k),\label{oct15-03}\\  
\langle e_{j_m},(\alpha,-1)\rangle=0&\Rightarrow& 
\langle e_{j_m},\alpha\rangle=0 \ \ \ (m=1,\ldots,s).\label{oct15-1-03} 
\end{eqnarray}
It is not hard to see that $v_{i_1},\ldots,v_{i_k},e_{j_1},\ldots,e_{j_s}$ 
are linearly independent vectors in $\mathbb{R}^n$. Indeed if 
$$
\lambda_1v_{i_1}+\cdots+\lambda_kv_{i_k}+\mu_1e_{j_1}+\cdots+\mu_se_{j_s}=0\
\ \ \   
(\lambda_{h},\mu_m\in\mathbb{R}),
$$
then taking inner product with $\alpha$ and using
Eqs.~(\ref{oct15-03}) and (\ref{oct15-1-03})  
we get
$$
\lambda_1+\cdots+\lambda_k=0\ \Rightarrow\ 
\lambda_1(v_{i_1},1)+\cdots+\lambda_k(v_{i_k},1)+\mu_1e_{j_1}+ 
\cdots+\mu_se_{j_s}=0. 
$$
Therefore $\lambda_{h}=0$ and $\mu_m=0$ for all ${h},m$, as desired.
From  
$\mathbb{R}_+{\cal A}'\subset H_{a_p}^+$ and
$H_{a_p}^+=H_{(\alpha,-1)}^+$ we get  
$\alpha\geq 0$ and $\langle\alpha,v_i\rangle\geq 1$ for all $i$.
Altogether  
we obtain that $\alpha$ is a basic feasible solution, that is,
$\alpha$  
is a vertex of $Q(A)$. \QED 

\medskip

Let $\mathfrak{p}_1,\ldots,\mathfrak{p}_s$ be the minimal primes 
of the edge ideal $I=I({\cal C})$ and let 
$$
C_k=\{x_i\vert\, x_i\in\mathfrak{p}_k\}\ \ \ \ (k=1,\ldots,s)
$$ 
be the corresponding minimal vertex covers of the clutter $\cal C$.
By  Proposition~\ref{1cover} and Theorem~\ref{gvvalencia-theo} in the
sequel we may assume that 

$$
a_k=(\textstyle\sum_{x_i\in C_k}e_i,-1)  \ \ \ \ (k=1,\ldots,s).
$$

\noindent{\it Notation\/} Let $d_k$ be the unique positive integer  
such that $d_ka_k$ has relatively prime integral entries. We set 
$\ell_k=d_ka_k$ for $k=1,\ldots,r$. If the 
first $n$ rational entries of $a_k$ are written in lowest terms, 
then $d_k$ is the least common multiple of the 
denominators. For $1\leq k\leq r$, 
we have $d_k=-\langle\ell_k,e_{n+1}\rangle$. 

\begin{Definition}\rm The set covering polyhedron $Q(A)$ is {\it
integral\/} 
if all its vertices 
have integral entries.  
\end{Definition}

\begin{Corollary}\label{1cover-1} The irreducible representation 
of the Rees cone has the form
\begin{equation}\label{okayama-oct} 
{\mathbb R}_+{\cal A}'=H_{e_1}^+\cap \cdots\cap H_{e_{n+1}}^+
\cap H_{\ell_1}^+\cap\cdots\cap H_{\ell_r}^+,
\end{equation}
$d_k=1$ if and only if $1\leq k\leq s$, and $Q(A)$ is integral if 
and only if $r=s$. 
\end{Corollary}

\demo It follows from Theorem~\ref{gvvalencia-theo} and
Corollary~\ref{may12-06}. \QED 

\medskip

\noindent {\it Notation} In the sequel we shall always assume 
that $\ell_1,\ldots,\ell_r$ are the integral vectors of
Eq.~(\ref{okayama-oct}).  

\medskip

Recall that the {\it Simis cone\/} of ${\cal A}$ is the rational 
polyhedral cone 
$$
{\rm Cn}({\cal A})=H_{e_1}^+\cap \cdots\cap H_{e_{n+1}}^+
\cap H_{\ell_1}^+\cap\cdots\cap H_{\ell_s}^+,
$$
and the {\it symbolic Rees algebra\/} of $I$ is the $K$ algebra:
$$
R_s(I)=R+I^{(1)}t+I^{(2)}t^2+\cdots+I^{(i)}t^i+\cdots\subset R[t],
$$
where $I^{(i)}=\mathfrak{p}_1^i\cap\cdots\cap\mathfrak{p}_s^i$ is the
$i${\it th\/} symbolic power of $I$. Symbolic Rees algebras have a
combinatorial interpretation \cite{cover-algebras}. Notice the
following 
description:
$$
I^{(b)}=(\{x^a\vert\, \langle (a,b),\ell_i\rangle\geq 0\mbox{ for
}i=1,\ldots,s\}).
$$

A first use of the Simis cone is the following expression 
for the symbolic Rees algebra. In particular $R_s(I)$ is a finitely 
generated $K$-algebra \cite{Lyu3} by Gordan's Lemma \cite{BHer}.

\begin{Theorem}[\cite{normali}]\label{oct24-03} 
If $S=\mathbb{Z}^{n+1}\cap{\rm Cn}({\cal A})$ and 
$K[S]=K[\{x^{a}t^b\vert\, (a,b)\in S\}]$ is its 
semigroup ring, then $R_s(I)=K[S]$.  
\end{Theorem}

Let $\mathbb{N}{\cal A}'$ be the subsemigroup of 
$\mathbb{N}^{n+1}$ generated by ${\cal A}'$, consisting of the 
linear combinations of ${\cal A}'$ with non-negative integer
coefficients. The Rees algebra of $I$ can be written as
\begin{eqnarray}
R[It]&=&
K[\{x^at^b\vert\, (a,b)\in\mathbb{N}{\cal A}'\}]\label{may6-06-2}\\ &=&
R\oplus It\oplus\cdots\oplus I^{i}t^i\oplus\cdots
\subset R[t].\label{may6-06-3}
\end{eqnarray}
According to \cite[Theorem~7.2.28]{monalg} and \cite[p.~168]{Vas1} the  
integral closure of $R[It]$ 
in its field of fractions can be expressed as
\begin{eqnarray}
\overline{R[It]}&=&K[\{x^at^b\vert\, (a,b)\in
\mathbb{Z}^{n+1}\cap \mathbb{R}_+{\cal A}'\}]\label{may6-06}\\ 
&=&R\oplus
\overline{I}t\oplus\cdots\oplus
\overline{I^i}t^i\oplus\cdots,\label{jun05-1-03} 
\end{eqnarray}
where $\overline{I^i}=(\{x^a\in R\vert\, \exists\, p\geq
1;(x^a)^{p}\in I^{pi}\})$ is the integral 
closure of $I^i$. Hence, by Eqs.~(\ref{may6-06-2}) to
$(\ref{jun05-1-03})$, we get that $R[It]$ is a normal domain if and
only if  
any of the following two equivalent conditions hold: 
\begin{description}
\item{\rm (a)} $\mathbb{N}{\cal A}'=
\mathbb{Z}^{n+1}\cap\mathbb{R}_+{\cal A}'$.\vspace{-1mm}
\item{\rm (b)} $I^{i}=\overline{I^i}$ for all $i\geq 1$.
\end{description}
If the second condition holds we say that $I$ is a {\it normal} 
ideal. 

\begin{Proposition} For $1\leq i\leq r$ we write $a_i=(a_i',-1)$. Let 
$B$ be the matrix with column vectors $a_1',\ldots,a_r'$ and let 
$Q=\mathbb{Q}_+^n+{\rm conv}(v_1,\ldots,v_q)$. Then 
\begin{description}
\item{\rm (a)} $\overline{I^i}=(\{x^a\in R\vert\,
a\in{iQ}\cap\mathbb{Z}^n\})$. 
\item{\rm (b)} $Q=Q(B)=\{x\vert\, x\geq 0;\,
xB\geq\mathbf{1}\}$. In particular $Q(B)$ is integral.
\end{description}
\end{Proposition}

\demo Part (a) follows from Eq.~(\ref{jun05-1-03}) and part (b) follows
from Eq.~(\ref{okayama}). \QED
 
\medskip

In the sequel $J_k^{(d_k)}$
will denote 
the ideal of $R[It]$ given by 
$$J_k^{(d_k)}=(\{x^at^b\in R[It]\vert\, \langle(a,b), 
\ell_k\rangle\geq d_k\}) \ \ \ \ (k=1,\ldots,r)$$
and $J_k$ will denote the ideal of $R[It]$ given by 
$$J_k=(\{x^at^b\in R[It]\vert\, \langle(a,b), 
\ell_k\rangle>0\}) \ \ \ \ (k=1,\ldots,r),$$
where $d_k=-\langle\ell_k,e_{n+1}\rangle$. If $d_k=1$, we have
$J_k^{(1)}=J_k$. In general $J_k^{(d_k)}$ might not be equal to the
$d_k${\it th} symbolic power of $J_k$. The localization of $R[It]$ at 
$R\setminus \mathfrak{p}_k$ is denoted by $R[It]_{\mathfrak{p}_k}$. 

\begin{Proposition}\label{oct28-03} $J_1,\ldots,J_r$ are height one
prime ideals containing $IR[It]$ and $J_k$ is equal to 
$\mathfrak{p}_kR[It]_{\mathfrak{p}_k}\cap R[It]$ for
$k=1,\ldots,s$. If $Q(A)$ is integral, then 
$$
{\rm rad}(IR[It])=J_1\cap J_2\cap\cdots\cap J_s.
$$ 
\end{Proposition}

\demo $IR[It]$ is clearly contained in $J_k$ for all $k$ by
construction. 
To show that $J_k$ is a prime ideal of height one it suffices to
notice that the  
right hand side of the isomorphism:
$$
R[It]/J_k\simeq K[\{x^at^b\in R[It]\vert\, \langle(a,b),\ell_k\rangle=0 \}]
$$
is an $n$-dimensional integral domain, because $F_k=\mathbb{R}_+{\cal
A}'\cap H_{\ell_k}$  
is a facet of the Rees cone for all $k$. Set 
$P_k=\mathfrak{p}_kR[It]_{\mathfrak{p}_k}\cap R[It]$ for $1\leq k\leq s$. 
This ideal 
is a minimal prime 
of $IR[It]$ (see \cite{HuSV}) and admits the following description 
\begin{eqnarray*}
P_k&=&
\mathfrak{p}_kR_{\mathfrak{p}_k}[\mathfrak{p}_kR_{\mathfrak{p}_k}t]
\cap R[It]
\\
&=&\mathfrak{p}_k+(\mathfrak{p}_k^2\cap I)t+
(\mathfrak{p}_k^3\cap I^2)t^2+
\cdots+(\mathfrak{p}_k^{i+1}\cap I^i)t^i+\cdots
\end{eqnarray*}
Notice that $x^a\in\mathfrak{p}_k^{b+1}$ if and only if 
$\langle a,\sum_{x_i\in C_k}e_i\rangle\geq b+1$. Hence $J_k=P_k$.

Assume that $Q(A)$ is integral, i.e., $r=s$. Take 
$x^\alpha t^b\in J_k$ for all $k$. Using Eq.~(\ref{okayama-oct}) it
is not hard  
to see that $(\alpha,b+1)\in\mathbb{R}_+{\cal A}'$, 
that is $x^\alpha t^{b+1}$  
is in $\overline{R[It]}$ and $x^\alpha t^{b+1}\in
\overline{I^{b+1}}t^{b+1}$. It follows  
that $x^\alpha t^b$ is a monomial in the radical of $IR[It]$. This proves 
the asserted equality.  \QED

\medskip

For use below recall that the {\it analytic spread\/} 
of $I$ is given by 
$$
\ell(I)=\dim R[It]/\mathfrak{m}R[It];\ \ \ \mathfrak{m}=(x_1,\ldots,x_n).
$$ 

\begin{Corollary}\label{march19-06} If $Q(A)$ is integral, then $\ell(I)<n$. 
\end{Corollary}

\demo Since $Q(A)$ is integral, we have $r=s$. If $\ell(I)=n$, 
then  the height of $\mathfrak{m}R[It]$ is equal to $1$. 
Hence there is a height one 
prime ideal $P$ of $R[It]$ such that 
$IR[It]\subset\mathfrak{m}R[It]\subset P$. 
By Proposition~\ref{oct28-03} the ideal $P$ has the form 
$\mathfrak{p}_kR[It]_{\mathfrak{p}_k}\cap R[It]$, this readily 
yields a contradiction. \QED 

\begin{Theorem}[\cite{Bur,E-H}]\label{burch} ${\rm inf}\{{\rm
depth}(R/I^i)\vert\, i\geq 1\}\leq\dim(R)-\ell(I)$. If\/ ${\rm
gr}_I(R)$ is 
Cohen-Macaulay, then the equality holds. \end{Theorem}

By a result of Brodmann \cite{brodmann}, the depth of $R/I^k$ 
is constant for $k$ sufficiently large. Broadmann improved
this inequality by showing that the constant value is bounded 
by $\dim(R)-\ell(I)$. For a study of the initial and limit behavior
of the numerical 
function $f(k)={\rm depth}\, R/I^k$ see \cite{depth}.

\begin{Lemma}\label{nov16-03} Let $x_1\in C_k$ for some $1\leq k\leq
s$. If $x^{v_i}=x_1x^{v_i'}$ for $1\leq i\leq p$ 
and $x_1\notin{\rm supp}(x^{v_i})$ for $i>p$, then there is 
$x^{v_j'}$ such that ${\rm supp}(x^{v_j'})\cap C_k=\emptyset$.  
\end{Lemma}

\demo If ${\rm supp}(x^{v_j'})\cap C_k\neq\emptyset$ for all $j$, 
then $C_k\setminus\{x_1\}$ is a vertex cover of $\cal C$, 
a contradiction because $C_k$ is a minimal vertex cover. \QED

\begin{Proposition} If $P\in\{J_1,\ldots,J_s\}$, then  
$R[It]\cap IR[It]_P=P$. 
\end{Proposition}

\demo Set $P=J_k$. We may assume that $x_1,\ldots,x_m$ (resp.
$x^{v_1}t,\ldots,x^{v_p}t$) is the set of all 
$x_i$ (resp. $x^{v_i}t$) such that $x_i\in P$ (resp. $x^{v_i}t\in P$).
Notice that $\mathfrak{p}_k$ is equal to $(x_1,\ldots,x_m)$ and set 
$C=\{x_1,\ldots,x_m\}$. In general the left hand side is contained
in $P$. To show the reverse  
inclusion we first prove the equality 
\begin{equation}\label{may2-06}
P=(x_1,\ldots,x_m,x^{v_1}t,\ldots,x^{v_p}t)R[It].
\end{equation}
Let $x^at^b\in P$. Thus
$x^at^b=x_1^{\mu_1}\cdots
x_n^{\mu_n}(x^{v_1}t)^{\lambda_1}\cdots(x^{v_q}t)^{\lambda_q}$ and 
$\langle (a,b),\ell_k\rangle>0$. Hence $\langle e_i,\ell_k\rangle>0$
for some $i$ or  
$\langle (v_j,1),\ell_k\rangle>0$ for some $j$. Therefore $x^at^b$
belongs 
to the right hand side of Eq.~(\ref{may2-06}), as required.

Case (I): Consider $x_\ell$ with $1\leq \ell\leq m$. By
Lemma~\ref{nov16-03} there is $j$  
such that $x^{v_j}=x_{\ell}x^{\alpha}$ and 
${\rm supp}(x^{\alpha})\cap C=\emptyset$. Thus since $x^{\alpha}$ is
not in  
$P$ (because of the second condition) we obtain $x_\ell\in R[It]\cap
IR[It]_P$. 

Case (II): Consider $x^{v_\ell}t$ with $1\leq \ell\leq p$. Since 
$$
\langle(v_\ell,1),e_1+\cdots+e_m-e_{n+1}\rangle\geq 1,
$$
the monomial $x^{v_\ell}$ contains at least two variables in $C$. Thus 
we may assume that $x_1,x_2$ are in the support of $x^{v_\ell}$. Again by 
Lemma~\ref{nov16-03} there are $j,j_1$ such that 
$x^{v_j}=x_1x^{\alpha}$, $x^{v_{j_1}}=x_2x^{\gamma}$, and the support of 
$x^\alpha$ and $x^\gamma$ disjoint from $C$. Hence the monomial 
$x^{v_\ell}x^{\alpha+\gamma}t$ belongs to $I^2t$ and
$x^{\alpha+\gamma}$ is not in $P$.  
Writing 
$$x^{v_\ell}t=({x^{v_\ell}x^{\alpha+\gamma}t})/{x^{\alpha+\gamma}},$$ 
we get $x^{v_\ell}t\in R[It]\cap IR[It]_P$. \QED

\begin{Lemma}
${\rm rad}(J_k^{(d_k)})=J_k$ for $1\leq k\leq r$.
\end{Lemma} 

\demo By construction one has ${\rm rad}(J_k^{(d_k)})\subset J_k$. The
reverse inclusion follows by noticing that if $x^at^b\in J_k$, 
then $(x^at^b)^{d_k}\in J_k^{(d_k)}$. \QED

\begin{Proposition}\label{primdec-irit} If $R[It]$ is normal, then 
$IR[It]=J_1^{(d_1)}\cap\cdots\cap J_r^{(d_r)}$.
\end{Proposition}

\demo ``$\subset$'': Let $x^a t^{b}\in IR[It]$. Since 
$x^a\in I^{b+1}$, we obtain $(a,b+1)\in\mathbb{N}{\cal A}'$. In
particular we get  
$(a,b+1)\in\mathbb{R}_+{\cal A}'$. Therefore
$$
0\leq \langle(a,b+1),\ell_k\rangle=\langle(a,b),\ell_k\rangle-d_k
$$ 
and consequently $x^a t^b\in J^{(d_k)}$ for $1\leq k\leq r$.

``$\supset$'': Let $x^a t^b\in J^{(d_k)}$ for all $k$. Since
$(a,b+1)\in \mathbb{R}_+{\cal A}'\cap\mathbb{Z}^{n+1}$, using that
$R[It]$ is normal yields  
$(a,b+1)\in\mathbb{N}{\cal A}'$. It follows that 
$x^a t^b\in I^{b+1}t^b\subset IR[It]$. \QED 

\medskip

A similar formula is shown in \cite{bruns-tania}.
The 
normality of $R[It]$ can be described in
terms of primary decompositions of $IR[It]$, see 
\cite[Proposition~2.1.3]{HSV}. 

\medskip

The following two nice 
formulas, pointed out to us by Vasconcelos, describe the difference between 
the symbolic Rees algebra of $I$ and the normalization of its Rees
algebra. If $\mathfrak{q}_k=J_k\cap{R}$ for $k=1,\ldots,r$, then 
\begin{eqnarray*}
R_s(I)=\bigcap_{k=1}^s R[It]_{\mathfrak{q}_k}\cap R[t];& &
\overline{R[It]}=\bigcap_{k=1}^r
\overline{R[It]}_{\mathfrak{q}_k}\cap R[t].
\end{eqnarray*}
These representations are linked to the so called Rees 
valuations of the ideal $I$, see \cite[Chapter~8]{bookthree} for further 
details.

\begin{Proposition}\label{oct30-03} The following conditions are
equivalent
\begin{description}
\item{\rm (a)} $Q(A)$ is integral.\vspace{-1mm}
\item{\rm (b)} ${\mathbb R}_+{\cal A}'=H_{e_1}^+\cap \cdots\cap H_{e_{n+1}}^+
\cap H_{\ell_1}^+\cap\cdots\cap H_{\ell_s}^+$, i.e., $r=s$.
\vspace{-1mm}
\item{\rm (c)} $R_s(I)=\overline{R[It]}$.\vspace{-1mm}
\item{\rm (d)} The minimal primes of $IR[It]$ are of the form 
$\mathfrak{p}_kR[It]_{\mathfrak{p}_k}\cap R[It]$.
\end{description}
\end{Proposition}

\demo (a) $\Leftrightarrow$ (b) $\Leftrightarrow$ (c): These
implications follow  
from Theorems~\ref{gvvalencia-theo} and \ref{oct24-03}. The other
implications  
follow readily using Proposition~\ref{oct28-03}. \QED

\begin{Definition}\rm Let $x^{u_k}=\prod_{x_i\in C_k}x_i$ for $1\leq
k\leq s$. The ideal  
of {\it vertex covers\/} of ${\cal C}$ is the ideal
$$
I_c({\cal C})=(x^{u_1},\ldots,x^{u_s})\subset R.
$$
The clutter of minimal vertex covers, denoted by $\cal D$ or 
$b({\cal C})$, is the {\it blocker\/} of $\cal C$.
\end{Definition}

In the literature $I_c({\cal C})$ is also called the {\it Alexander
dual\/} of 
$I$ because if $\Delta$ is the Stanley-Reisner complex of $I$, then
$I_c({\cal C})$ is the Stanley-Reisner ideal 
of the Alexander dual of $\Delta$. The survey 
article \cite{herzog-survey} explains the role of Alexander duality 
to prove combinatorial and algebraic theorems.  

\begin{Example}\rm Let $I=(x_1x_2x_5,x_1x_3x_4,x_2x_3x_6,x_4x_5x_6)$.
The clutter of $I$ is denoted by ${\cal Q}_6$. Using {\it Normaliz\/}
\cite{B}  
and  Proposition~\ref{oct30-03} we obtain:  
$$R[It]\subsetneq{R_s(I)}=\overline{R[It]}=R[It][x_1\cdots x_6t^2]\
\mbox{ and }\ R[I_c({\cal Q}_6)t]=\overline{R[I_c({\cal Q}_6)t]}.$$
\end{Example}

\begin{Proposition}[\cite{alexdual}]\label{nov3-03}
If $\overline{R[It]}=R_s(I)$ and $J=I_c({\cal C})$, then
$\overline{R[Jt]}=R_s(J)$. 
\end{Proposition}

\begin{Corollary}{\rm \cite[Theorem~1.17]{cornu-book}} If $Q(A)$ is
integral and $A'$ is the incidence matrix of the clutter of minimal 
vertex covers of $\cal C$, then $Q(A')$ is integral. 
\end{Corollary}

\demo It follows at once from Propositions~\ref{oct30-03} and
\ref{nov3-03}. \QED

\begin{Definition}\rm Let 
$X'=\{x_{i_1},\ldots,x_{i_r},x_{j_1},\ldots,x_{j_s}\}$ be a subset of $X$. 
A {\it minor\/} of $I$ is a proper ideal $I'$ of $R'=K[X\setminus
X']$ obtained from  
$I$ by making 
$x_{i_k}=0$ and $x_{j_\ell}=1$ for all $k,\ell$. The ideal $I$ is
considered itself a minor.  
A {\it minor\/} of $\cal C$ 
is a clutter ${\cal C}'$ that corresponds to a minor $(0)\subsetneq
I'\subsetneq R'$.  
\end{Definition}

Notice that ${\cal C}'$ is obtained from $I'$ by considering the unique set 
of square-free monomials of $R'$ that minimally generate $I'$.

\begin{Proposition} If $\overline{I^i}=I^{(i)}$ for some $i\geq 2$
and $J=I'$ is a  
minor of $I$, then $\overline{J^i}=J^{(i)}$.
\end{Proposition}

\demo Assume that $J$ is the minor obtained from $I$ by making $x_1=0$. 
Take $x^a\in{J^{(i)}}$. Then $x^a\in I^{(i)}=\overline{I^i}$ because 
$J\subset I$. Thus $x^a\in\overline{I^i}$. Since 
$x_1\notin{\rm supp}(x^a)$ it follows that $x^a\in \overline{J^i}$. This 
proves $J^{(i)}\subset\overline{J^i}$. The other inclusion is clear because 
$J^{(i)}$ is integrally closed. 

Assume that $J$ is the minor obtained from $I$ by making $x_1=1$. Take 
$x^a\in J^{(i)}$. Notice that $x_1^ix^a\in I^{(i)}=\overline{I^i}$. 
Indeed if 
$x_1\in\mathfrak{p}_k$, then $x_1^i\in\mathfrak{p}_k^i$, and if 
$x_1\notin\mathfrak{p}_k$, then $J\subset \mathfrak{p}_k$ and 
$x^a\in\mathfrak{p}_k^i$. Since 
$x_1\notin{\rm supp}(x^a)$ it follows that $x^a\in \overline{J^i}$.\QED

\begin{Corollary}\label{nov19-1-03} If $R_s(I)=\overline{R[It]}$, 
then $R_s(I')=\overline{R'[I't]}$ 
for any minor $I'$ of $I$.
\end{Corollary}

\begin{Proposition}\label{nov2-03} Let ${\cal D}$ be the clutter of
minimal vertex  
covers of $\cal C$. If 
$\overline{R[It]}$ is equal to $R_s(I)$ 
and $|A\cap B|\leq 2$ for $A\in \cal C$ and $B\in \cal D$, 
then $R[It]$ is normal. 
\end{Proposition}

\demo Let $x^at^b=x_1^{a_1}\cdots x_n^{a_n}t^b\in\overline{R[It]}$ be
a minimal generator,  
that is $(a,b)$ cannot be written as a sum of two non-zero integral
vectors in the Rees cone  
$\mathbb{R}_+{\cal A}'$. We may assume $a_i\geq 1$ 
for $1\leq i\leq m$, $a_i=0$ for $i>m$, and $b\geq 1$. 

Case (I): $\langle(a,b),\ell_i\rangle>0$ for all $i$. The vector
$\gamma=(a,b)-e_1$  
satisfies $\langle\gamma,\ell_i\rangle\geq 0$ for all $i$, that is 
$\gamma\in\mathbb{R}_+{\cal A}'$. Thus since $(a,b)=e_1+\gamma$ we derive a 
contradiction. 

Case (II): $\langle(a,b),\ell_i\rangle=0$ for some $i$. We may assume
$$
\{\ell_i\vert\, \langle(a,b),\ell_i\rangle=0\}=\{\ell_1,\ldots,\ell_p\}.
$$

Subcase (II.a): $e_i\in H_{\ell_1}\cap\cdots\cap H_{\ell_p}$ for some
$1\leq i\leq m$. 
It is not hard to verify that the vector 
$\gamma=(a,b)-e_i$ satisfies $\langle\gamma,\ell_k\rangle\geq 0$ for all 
$1\leq k\leq s$. Thus $\gamma\in\mathbb{R}_+{\cal A}'$, a
contradiction because  
$(a,b)=e_i+\gamma$. 

Subcase (II.b): $e_i\notin H_{\ell_1}\cap\cdots\cap H_{\ell_p}$ for
all $1\leq i\leq m$. Since 
the vector $(a,b)$ belongs to $\mathbb{R}_+{\cal A}'$, it follows
(see the proof of  
Theorem~\ref{oct31-03}) that we can write
\begin{equation}\label{may6-06-1}
(a,b)=\lambda_1(v_1,1)+\cdots+\lambda_q(v_q,1) \ \ \ \ (\lambda_i\geq 0).
\end{equation}
By the choice of $x^at^b$ we may assume $0<\lambda_1<1$. Set 
$\gamma=(a,b)-(v_1,1)$ and notice that by Eq.~(\ref{may6-06-1}) this vector has 
non-negative entries. We claim that $\gamma$ is in the Rees cone. Since 
by hypothesis one has $0\leq \langle(v_1,1),\ell_j\rangle\leq 1$ for all $j$ we 
readily obtain
$$
\langle\gamma,\ell_k\rangle=\left\{\begin{array}{cccl}
\langle(a,b),\ell_k\rangle-\langle (v_1,1),\ell_k\rangle&=& 0
&\mbox{if }\ 1\leq k\leq p,\\ 
\langle(a,b),\ell_k\rangle-\langle(v_1,1),\ell_k\rangle&\geq& 0 &\mbox{otherwise.} 
\end{array} \right.
$$
Thus $\gamma\in\mathbb{R}_+{\cal A}'$ and $(a,b)=(v_1,1)+\gamma$. As
a result $\gamma=0$  
and $x^at^b\in R[It]$, as desired. \QED

\section{K\"onig property of clutters and normality}

Let us introduce a little bit more notation and definitions. Recall
that the  
{\it Ehrhart 
ring\/} of the lattice polytope $P={\rm conv}(v_1,\ldots,v_q)$ is the
subring  
$$
A(P)=K[\{x^at^i\vert\, a\in \mathbb{Z}^n \cap iP; i\in \mathbb{N}\}]\subset R[t],
$$
and the homogeneous {\it monomial subring\/} generated by 
$Ft=\{x^{v_1}t,\ldots,x^{v_q}t\}$ over the field $K$ is the subring
$K[Ft]\subset R[t]$.  

\begin{Theorem}\label{oct31-03} If $\overline{R[It]}=R_s(I)$ and 
$K[Ft]=A(P)$, then $R[It]$ is normal.
\end{Theorem}

\demo Let $x^at^b=x_1^{a_1}\cdots x_n^{a_n}t^b\in\overline{R[It]}$ be
a minimal generator,  
that is $x^at^b$ cannot be written as a product of 
two non-constant monomials of $\overline{R[It]}$. We may assume $a_i\geq 1$ 
for $1\leq i\leq m$, $a_i=0$ for $i>m$, and $b\geq 1$.

Case (I): $\langle(a,b),\ell_i\rangle>0$ for all $i$. The vector $\gamma=(a,b)-e_1$ 
satisfies $\langle\gamma,\ell_i\rangle\geq 0$ for all $i$, that is 
$\gamma\in\mathbb{R}_+{\cal A}'$. Thus since $x_1$ and 
$x_1^{a_1-1}x_2^{a_2}\cdots x_{n}^{a_n}t^b$ are in $\overline{R[It]}$ we get a 
contradiction. In conclusion this case cannot occur. 

Case (II): $\langle(a,b),\ell_i\rangle=0$ for some $i$. We may assume
$$
\{\ell_i\vert\, \langle(a,b),\ell_i\rangle=0\}=\{\ell_1,\ldots,\ell_p\}.
$$

Subcase (II.a): $e_i\in H_{\ell_1}\cap\cdots\cap H_{\ell_p}$ for some
$1\leq i\leq m$. For  
simplicity of notation assume $i=1$. The vector $\gamma=(a,b)-e_1$ satisfies
$$
\langle\gamma,\ell_k\rangle=\left\{\begin{array}{cccl}
\langle(a,b),\ell_k\rangle-\langle e_1,\ell_k\rangle&=& 0 &\mbox{if
}\ 1\leq k\leq p,\\ 
\langle(a,b),\ell_k\rangle-\langle e_1,\ell_k\rangle&\geq& 0 &\mbox{otherwise.} 
\end{array} \right.
$$
Thus $\gamma\in\mathbb{R}_+{\cal A}'$. Proceeding as in Case (I) we
derive a contradiction. 

Subcase (II.b): $e_i\notin H_{\ell_1}\cap\cdots\cap H_{\ell_p}$ for
all $1\leq i\leq m$. The  
vector $(a,b)$ belongs to the polyhedral cone
$$
C=H_{\ell_1}\cap\cdots\cap H_{\ell_p}\cap \mathbb{R}_+{\cal A}'.
$$
Hence we can write
\begin{eqnarray*}
(a,b)&=&\lambda_1(v_1,1)+\cdots+\lambda_q(v_q,1)+\mu_1e_1+\cdots+\mu_ne_n
\ \ \ \  
(\lambda_i;\mu_j\geq 0),\\ 
\langle(a,b),\ell_k\rangle&=&\lambda_1\langle(v_1,1),\ell_k\rangle+\cdots
+\lambda_q\langle(v_q,1),\ell_k\rangle+\\ 
& &\mu_1\langle e_1,\ell_k\rangle+\cdots+
\mu_n\langle e_n,\ell_k\rangle=0
\end{eqnarray*}
for $k=1,\ldots,p$. From the first equality we get $\mu_i=0$ for
$i>m$ because  
$a_i=0$ for $i>m$. If $\mu_i>0$ for some $1\leq i\leq m$, then
$\langle e_i,\ell_k\rangle=0$ for  
$1\leq k\leq p$, a contradiction. Hence $\mu_i=0$ for all $i$.
Therefore $a/b\in P$ and  
$a\in\mathbb{Z}^{n}\cap bP$. This proves $x^at^b\in A(P)=K[Ft]\subset
R[It]$, as  
desired. \QED 

\begin{Proposition}[\cite{matrof}] If $\ell_k=d_ka_k$ has the form 
$$
\ell_k=e_{i_1}+\cdots+e_{i_k}-d_ke_{n+1}\ \ \ \ \ (1\leq
i_1<\cdots<i_k\leq n)
$$
for $k=1,\ldots,r$, then $A(P)[x_1,\ldots,x_n]=\overline{R[It]}$.
\end{Proposition}

\demo The Ehrhart ring is contained in $\overline{R[It]}$. Thus the equality 
follows using the proof of Theorem~\ref{oct31-03}. \QED

\begin{Proposition} If $x^{v_1},\ldots,x^{v_q}$ have degree $d\geq 2$ and 
$\overline{I^b}=I^{(b)}$ for all $b$, then $\overline{I^b}$ is generated by 
monomials of degree $bd$ for $b\geq 1$. 
\end{Proposition}

\demo The monomial ideal $\overline{I^b}$ has a unique minimal set of
generators consisting  
of monomials. Take $x^a$ in this minimal set. Notice that 
$(a,b)\in\mathbb{R}_+{\cal A}'$. Thus we may proceed as in the proof of 
Theorem~\ref{oct31-03} to obtain that $(a,b)$ is in the cone generated by 
$\{(v_1,1),\ldots,(v_q,1)\}$. This yields $\deg(x^a)=bd$. \QED

\begin{Proposition}\label{nov1-03} If $x^{v_1},\ldots,x^{v_q}$ have
degree $d\geq 2$,  
then $I^{i}=I^{(i)}$ for all $i\geq 1$ if and only if $Q(A)$ is integral and 
$K[Ft]=A(P)$.
\end{Proposition}

\demo $\Rightarrow$) By Proposition~\ref{oct30-03} the polyhedron
$Q(A)$ is integral.  
Since $I^{(i)}$ is integrally closed \cite[Corollary~7.3.15]{monalg},
we get that  
$R[It]$ is normal. Therefore applying \cite[Theorem~3.15]{ehrhart} we obtain 
$K[Ft]=A(P)$, here the hypothesis on the degrees of $x^{v_i}$ is essential.

$\Leftarrow$) By Proposition~\ref{oct30-03}
$\overline{I^{i}}=I^{(i)}$ for all $i$,  
thus applying Theorem~\ref{oct31-03} gives $R[It]$ normal and we get
the  required equality. Here the hypothesis on the degrees 
of $x^{v_i}$ is not needed. \QED

\begin{Definition}\rm The clutter $\cal C$ satisfies the {\it max-flow min-cut\/}
(MFMC)
property if both sides 
of the LP-duality equation
\begin{equation}\label{jun6-2-03-1}
{\rm min}\{\langle \alpha,x\rangle \vert\, x\geq 0; xA\geq \mathbf{1}\}=
{\rm max}\{\langle y,\mathbf{1}\rangle \vert\, y\geq 0; Ay\leq\alpha\} 
\end{equation}
have integral optimum solutions $x$ and $y$ for each non-negative
integral vector $\alpha$.  
\end{Definition}

It follows from \cite[pp.~311-312]{Schr} that $\cal C$ has the
MFMC property if and only 
if the maximum in Eq.~(\ref{jun6-2-03-1}) has an optimal integral
solution $y$ for each  
non-negative integral vector $\alpha$. Thus the system $x\geq 0;\
xA\geq\mathbf{1}$ is TDI  if and only if $\cal C$ has the max-flow
min-cut property.

\medskip

A ring is called 
{\it reduced} if $0$ is its only nilpotent element. For convenience
let us state some  
known characterizations of the reducedness of the associated graded ring.

\begin{Theorem}[\cite{normali,clutters,HuSV}]\label{noclu} The
following conditions  
are equivalent
\begin{description}
\item{\rm(i)\ \ } ${\rm gr}_I(R)$ is reduced.
\vspace{-1mm}
\item{\rm (ii)\ } $R[It]$ is normal and $Q(A)$ is an integral
polyhedron.\vspace{-1mm} 
\item{\rm (iii)} $I$ is normally torsion free, that is, $I^i=I^{(i)}$ 
for all $i\geq 1$.\vspace{-1mm}
\item{\rm (iv)} $x\geq 0;\, xA\geq \mathbf{1}$ is a {\rm TDI}
system.\vspace{-1mm} 
\item{\rm (iv)} $\cal C$ has the max-flow min-cut property. 
\end{description}
\end{Theorem}

\begin{Corollary}{\rm \cite[Theorem~1.3]{CGM}} Let ${\cal D}$ be the
clutter of minimal vertex  
covers of $\cal C$. If $Q(A)$ is integral and 
$|A\cap B|\leq 2$ for $A\in \cal C$ and $B\in \cal D$, then 
$x\geq 0;\, xA\geq \mathbf{1}$ is a {\rm TDI} system. 
\end{Corollary}

\demo By Proposition~\ref{nov2-03} the Rees algebra $R[It]$ is normal.
To complete the proof apply Theorem~\ref{noclu}. \QED

\begin{Lemma}\label{dec25-02}
If $I$ is a monomial ideal of $R$, then the nilradical of
the associated graded ring of $I$ is given by
$$
{\rm nil}({\rm gr}_I(R))
=(\{\overline{x^\alpha}\in I^i/I^{i+1}\vert\, 
x^{s\alpha}\in I^{si+1}; i\geq 0; s\geq 1\}).
$$
\end{Lemma}

\demo The nilradical of ${\rm gr}_I(R)$ is graded with respect to 
the fine grading, and thus it is generated by homogeneous 
elements. \QED 

\begin{Definition}\rm The matrix $A$ is {\it balanced\/} if $A$ has no
square submatrix of odd order with exactly two $1$'s in
each row and column. $A$ is {\it totally unimodular\/} if each
$i\times i$ minor of $A$ is $0$ or $\pm 1$ for all $i\geq 1$.
\end{Definition}

\begin{Proposition}\label{dec30-03} If $A$ is balanced, 
then ${\rm gr}_{I}(R)$ is 
reduced.  
\end{Proposition}

\demo Let $\overline{x^{\alpha}}\in I^i/I^{i+1}$ be in ${\rm
nil}({\rm gr}_I(R))$,  that is  $x^{s\alpha}\in I^{is+1}$ for 
some $0\neq s\in\mathbb{N}$. By Lemma~\ref{dec25-02} we need only show
 $\overline{x^\alpha}=\overline{0}$. It follows rapidly that the
maximum in Eq.~(\ref{jun6-2-03}) 
is greater or equal than $i+1/s$. By 
\cite[Theorem~21.8, p.~305]{Schr} the maximum in
Eq.~(\ref{jun6-2-03}) 
has an integral optimum solution $y=(y_1,\ldots,y_q)$. Thus 
$y_1+\cdots+y_q\geq i+1$. Since $y$ satisfies $y\geq 0$ and $Ay\leq
\alpha$ 
we obtain $x^\alpha\in I^{i+1}$. This proves
$\overline{x^\alpha}=\overline{0}$, as required. \QED

\begin{Proposition}\label{balanced-mfmc-blocker} If $A$ is balanced 
and $J=I_c({\cal C})$, then $R[Jt]=R_s(J)$. 
\end{Proposition}

\demo Let $\cal D$ be the blocker of $\cal C$. By 
\cite[Corollary~83.1a(v), p.~1441]{Schr2}, we get that $\cal D$
satisfy the max-flow min-cut
property. Thus the equality follows at 
once from Theorem~\ref{noclu}. \QED

\begin{Proposition}\label{oct20-03} If ${\rm gr}_I(R)$ is reduced 
$($resp. $R[It]$ is normal$)$ and $I'$ is a 
minor of $I$, then ${\rm gr}_{I'}(R')$ is reduced $($resp. $R'[I't]$
is normal$)$. 
\end{Proposition}

\demo Notice that we need only show the result when $I'$ is a minor 
obtained from $I$ by making $x_1=0$ or $x_1=1$. Using
Lemma~\ref{dec25-02} both  
cases are quite easy to prove. \QED

\begin{Definition}\rm A clutter $\cal C$ satisfies the {\it packing
property\/}  
(PP) if all its minors satisfy the K\"onig property, 
that is, $\alpha_0({\cal C}')=\beta_1({\cal C}')$ 
for any minor ${\cal C}'$ of $\cal C$.
\end{Definition}

\begin{Corollary}\label{oct21-03} If the ring ${\rm gr}_I(R)$ is 
reduced, then $\alpha_0({\cal C}')=
\beta_1({\cal C}')$ for any minor ${\cal C}'$ of $\cal C$.
\end{Corollary}

\demo Let ${\cal C}'$ be any minor of $\cal C$ and let $I'$ be its
clutter ideal. We denote the incidence matrix of ${\cal C}'$ by $A'$.
By Proposition~\ref{oct20-03} the 
associated graded ring ${\rm gr}_{I'}(R')$ 
is reduced. Hence by Theorem~\ref{noclu} the clutter ${\cal C}'$ has 
the max-flow min-cut property. In particular the LP-duality equation
$$
{\rm min}\{\langle\mathbf{1},x\rangle \vert\, x\geq 0; xA'\geq \mathbf{1}\}=
{\rm max}\{\langle y,\mathbf{1}\rangle \vert\, y\geq 0; A'y\leq\mathbf{1}\}
$$
has optimum integral solutions $x$, $y$. To complete the 
proof notice that the left hand side of this 
equality is $\alpha_0({\cal C}')$ and the right hand side 
is $\beta_1({\cal C}')$. \QED

\begin{Remark}\rm If $I$ is the facet ideal of a simplicial tree, 
then ${\rm gr}_I(R)$ is reduced. This follows from \cite[p.~174]{faridi} 
using the proof of \cite[Corollary~3.2, p.~399]{ITG}. In particular
$\cal C$ has the  
K\"onig property, this was shown in \cite[Theorem~5.3]{faridi1}.  
\end{Remark}

\begin{Corollary}\label{jan20-05} If $\cal C$ has the max-flow
min-cut property,  
then $\cal C$ has the packing property. 
\end{Corollary}

\demo It follows at once from Theorem~\ref{noclu} and 
Corollary~\ref{oct21-03}. \QED

\medskip

Conforti and Cornu\'ejols conjecture that the converse is also
true: 

\begin{Conjecture}{\rm \cite[Conjecture~1.6]{cornu-book}}
\label{conforti-cornuejols1}\rm 
If the clutter $\cal C$ has the
packing property, then $\cal C$ has the max-flow min-cut property.
\end{Conjecture}

Next we state the converse of Corollary~\ref{oct21-03} as an
algebraic version of 
this interesting conjecture which to our best knowledge 
is still open:  

\begin{Conjecture}\label{conforti-cornuejols}\rm If 
$\alpha_0({\cal C}')=\beta_1({\cal C}')$ for 
all minors $\cal C'$ of $\cal C$, then the 
ring ${\rm gr}_I(R)$ is reduced.
\end{Conjecture} 

It is known \cite[Theorem~1.8]{cornu-book} that clutters with the packing 
property have integral set covering polyhedrons. As a consequence,
using Theorem~\ref{noclu},  
this conjecture reduces to the following:  

\begin{Conjecture}\label{con-cor-vila}\rm If $\alpha_0({\cal
C}')=\beta_1({\cal C}')$ for  
all minors $\cal C'$ of $\cal C$, then $R[It]$ is normal.
\end{Conjecture} 

In this paper we will give some support for this conjecture using an 
algebraic approach. 

\begin{Proposition}\label{march18-06-1} 
Let $J_i$ be the ideal obtained from $I$ by making $x_i=1$. If $Q(A)$
is an integral polyhedron, then the ideal $I$ is normal if and only
if $J_i$ 
is normal for all
$i$ and ${\rm depth}(R/I^k)\geq 1$ for all $k\geq 1$.
\end{Proposition}

\demo $\Rightarrow$) The normality of an ideal is closed under minors
\cite[Proposition~4.3]{normali}, hence $J_i$ is normal for all $i$. Using
Theorem~\ref{burch} and Corollary~\ref{march19-06}  we get
that ${\rm depth}(R/I^i)\geq 1$ for all $i$. 

$\Leftarrow$) It follows readily by adapting the arguments given in the proof 
of the normality criterion \cite[Theorem~4.4]{normali}.\QED

\medskip

By Proposition~\ref{march18-06-1} we obtain that
Conjecture~\ref{conforti-cornuejols} also reduces to: 

\begin{Conjecture}\rm If $\alpha_0({\cal C}')=\beta_1({\cal C}')$ for
any minor ${\cal C}'$ of $\cal C$, then
$$
{\rm depth}(R/I^i)\geq 1\ \mbox{ for all }i\geq 1.
$$
\end{Conjecture}

\noindent {\it Notation} For an integral matrix $B\neq(0)$, the 
greatest common divisor of all the nonzero $r\times r$
subdeterminants of $B$ will be  
denoted by $\Delta_r(B)$.

\begin{Corollary} If $x^{v_1},\ldots,x^{v_q}$ are monomials of degree
$d\geq 2$  
such that ${\rm gr}_I(R)$ is reduced and the matrix  
$$
B=\left(\hspace{-1mm}
\begin{array}{ccc}
v_1&\cdots&v_q\\ 
1&\cdots &1
\end{array}\hspace{-1mm}
\right)
$$
has rank $r$, then $\Delta_r(B)=1$ and $B$ diagonalizes over
$\mathbb{Z}$ to an identity
matrix.
\end{Corollary}

\demo By Proposition~\ref{nov1-03} we obtain $A(P)=K[Ft]$. Hence a direct 
application of \cite[Theorem~3.9]{ehrhart} gives $\Delta_r(B)=1$. \QED

\medskip

This result suggest the following weaker conjecture of Villarreal:  

\begin{Conjecture} If $\alpha_0({\cal C}')=\beta_1({\cal C}')$ for 
all minors $\cal C'$ of $\cal C$ and $x^{v_1},\dots,x^{v_q}$ 
have degree $d\geq 2$, then $\Delta_r(B)=1$, where $r={\rm rank}(B)$. 
\end{Conjecture}

Let $G$ be a matroid on $X$ of rank $d$ and let $\cal B$ be the
collection of 
bases of $G$. The set of all 
square-free monomials $x_{i_1}\cdots x_{i_d}\in R$ such 
that $\{x_{i_1},\ldots,x_{i_d}\}\in {\cal B}$ will be denoted by
$F_G$ and the subsemigroup  
(of the multiplicative semigroup of monomials of $R$) generated by
$F_G$ will be denoted  
by $\mathbb{M}_G$. The {\it basis monomial ring\/} of $G$ is the
monomial subring $K[F_G]=K[\mathbb{M}_G]$. The ideal $I({\cal
B})=(F_G)$ is called the  
{\it basis monomial ideal} of $G$. 
An open problem in the area is whether the toric
ideal of $K[F_G]$ is generated by quadrics, see 
\cite[Conjecture~12]{white1}. 
This has been shown for graphic matroids \cite{blasiak}.

\medskip

The next result implies the normality of the basis monomial ring of
$G$.  
\begin{Proposition}[\cite{white}]\label{nwhite} If $x^a$ is a
monomial of degree $\ell{d}$ for some $\ell\in\mathbb{N}$ such that
$(x^{a})^p\in\mathbb{M}_G$ for some  
$0\neq p\in\mathbb{N}$, then $x^a\in\mathbb{M}_G$.
\end{Proposition}

\begin{Proposition} If $I=I({\cal B})$ and $\cal B$ 
satisfies the packing property, then ${\rm gr}_I(R)$ is reduced.
\end{Proposition}

\demo First we show the equality $A(P)=K[F_Gt]$. It suffices to prove
the inclusion $A(P)\subset K[F_Gt]$. Take 
$x^at^b\in A(P)$, i.e., $x^a\in\mathbb{Z}^n\cap bP$. Hence $x^a$ has 
degree $bd$ and $(x^{a})^p\in\mathbb{M}_G$ for some positive integer
$p$. By  
Proposition~\ref{nwhite} we get $x^a\in\mathbb{M}_G$. It is seen that 
$x^at^b$ is in $K[F_Gt]$. Since $Q(A)$ is integral
\cite{cornu-book}, using  
Theorem~\ref{oct31-03} we get that $R[It]$ is normal. Thus both
conditions yield  
that ${\rm gr}_I(R)$ is reduced according to Theorem~\ref{noclu}. \QED 

\medskip

This proof can be simplified using that the basis monomial ideal of a
matroid is normal \cite{matrof}.

\begin{Corollary} Let $X_1,\ldots,X_d$ be a family of disjoint
sets of variables and let $M$ be the transversal
matroid whose collection of basis is 
$$
{\cal C}=\{\{y_1,\ldots,y_d\}\vert\, y_i\in X_i\, \forall\, i\}.
$$ 
If $I=I({\cal C})$, then ${\rm gr}_I(R)$ is reduced. 
\end{Corollary}

The combinatorial equivalencies in the next result are well known
\cite{cornu-book,CGM}.  Our contribution here is to link the reducedness of
the associated graded ring with the integrality of $Q(A)$. 

\begin{Proposition} If $\cal C$ is a simple graph, then the following 
are equivalent{\rm :}
\begin{description}
\item{\rm (a)} ${\rm gr}_I(R)$ is reduced.\vspace{-1mm}
\item{\rm (b)} $\cal C$ is bipartite.\vspace{-1mm}
\item{\rm (c)} $Q(A)$ is integral.\vspace{-1mm}
\item{\rm (d)} $\cal C$ has the packing property. 
\end{description}
\end{Proposition}

\demo By \cite[Theorem~5.9]{ITG} (a) and (b) are equivalent. Applying 
Theorem~\ref{noclu} and Proposition~\ref{nov2-03} we obtain that (a) is 
equivalent to (c). By Corollary~\ref{oct21-03} condition (a) implies
(d). Finally using  
\cite[Theorem~1.8]{cornu-book} we obtain that (d) implies (c). \QED

\begin{Corollary}[\cite{alexdual}]\label{nov9-2-03} If $\cal C$ is a
bipartite graph and $J=I_c({\cal C})$,  
then ${\rm gr}_J(R)$ is reduced. 
\end{Corollary}

\demo The matrix $A$ is totally unimodular
\cite[p.~273]{Schr}, hence 
$Q(A)$ is integral. By Proposition~\ref{nov3-03} we get 
$\overline{R[Jt]}=R_s(J)$. On the other hand $R[Jt]$ is normal by 
Proposition~\ref{nov2-03}. Thus by Theorem~\ref{noclu} the ring ${\rm
gr}_J(R)$ is reduced. \QED

\begin{Proposition}
If $\cal C$ has the packing property and $I=I({\cal C})$, then
$I^2=\overline{I^2}$.  
\end{Proposition}

\demo By induction on $n$. Assume $\overline{I^2}\neq I^2$ and
consider  
$M=\overline{I^2}/I^2$. If
$\mathfrak{p}\neq\mathfrak{m}=(x_1,\ldots,x_n)$ is  
a prime ideal of $R$, then by induction $M_\mathfrak{p}=(0)$. Thus
$\mathfrak{m}$  
is the only associated prime of $M$ and there is an embedding 
$R/{\mathfrak m}\hookrightarrow M$, $\overline{1}\mapsto
\overline{x^{a}}$, where $\overline{x^{a}}\in\overline{I^2}\setminus
I^2$ and $x_ix^{{a}}\in I^2$  
for all $i$. Notice that by induction all the entries ${a}_i$ of
${a}$ are  
positive. We consider two cases. Assume $a_i\geq 2$ for some $i$, say
$i=1$. Given a monomial  
$x^\alpha$, the monomial obtained from $x^\alpha$ 
by making $x_1=1$ is denoted by $x^{\alpha'}$. Then making $x_1=1$
and using that $x_1x^a\in I^2$ gives $x^{a'}=x^{v_1'}x^{v_2'}x^\delta$, 
hence $x^a=x_1^{a_1}x^{a'}=x^{v_1}x^{v_2}x^\gamma\in I^2$, a
contradiction. On the other hand  
if $a_i=1$ for all $i$, then $x^a=x_1\cdots x_n\in I^g\subset I^2$, where 
$g={\rm ht}(I)$, a contradiction. Therefore $I^2=\overline{I^2}$. \QED

\medskip

Recall that $I$ is said to be {\it unmixed\/} if all the minimal
vertex covers  
of $\cal C$ 
have the same cardinality.

\begin{Lemma} If $I$ is an unmixed ideal and $\cal C$ satisfies the
K\"onig property,  
then $x^{\mathbf 1}=x_1x_2\cdots x_n$ belongs to the subring
$K[x^{v_1},\ldots,x^{v_q}]$.  
\end{Lemma}

\demo We may assume $x^{\mathbf 1}=x^{v_1}\cdots x^{v_g}x^\delta$,
where $g$ is the  
height of $I$. If $\delta\neq 0$, pick $x_n\in{\rm supp}(x^\delta)$.
Since the  
variable $x_n$ occurs in some monomial of $I$, there is a minimal
prime $\mathfrak{p}$ containing $x_n$. Thus  
using that $x^{v_1},\ldots,x^{v_g}$ have disjoint supports we
conclude that $\mathfrak{p}$ 
 contains at least $g+1$ variables, a contradiction. \QED

\begin{Proposition}\label{nov20-1-03} Let $I_i=I\cap
K[X\setminus\{x_i\}]$. If $I$  
is an unmixed ideal such that the following conditions hold
\begin{description}
\item{\rm ($\mathrm{a}_1$)} $Q(A)$ is integral,\vspace{-1mm}
\item{\rm ($\mathrm{a}_2$)} $I_i$ is normal for $i=1,\ldots,n$,
and\vspace{-1mm} 
\item{\rm ($\mathrm{a}_3$)} $\cal C$ has the K\"onig property,  
\end{description}
then $R[It]$ is normal.
\end{Proposition}

\demo Take $x^at^b=x_1^{a_1}\cdots x_n^{a_n}t^b\in\overline{R[It]}$ a
minimal generator. By the second condition we may assume 
$a_i\geq 1$ for all $i$. Set $g={\rm ht}(I)$. Notice that $x_1\cdots
x_nt^g$ is  
in $\overline{R[It]}$ because $Q(A)$ is integral, 
this follows from Corollary~\ref{1cover-1} and Eq.~(\ref{may6-06}). 
We claim that $b\leq g$. If $b>g$, consider the decomposition
$$
x^at^b=(x_1\cdots x_nt^g)(x_1^{a_1-1}\cdots x_n^{a_n-1}t^{b-g}).
$$  
To derive a contradiction consider the irreducible representation 
of the Rees cone $\mathbb{R}_+{\cal A}'$ given in
Eq.~(\ref{okayama-oct}).  
Observe that 
$$
\textstyle\sum_{x_i\in C_k}a_i\geq b\ \ \ \ \ (k=1,\ldots,s)
$$
because $(a,b)\in\mathbb{R}_+{\cal A}'$. Now since $I$ 
is unmixed we get
$$
\textstyle\sum_{x_i\in C_k}(a_i-1)\geq b-g\ \ \ \ \ (k=1,\ldots,s),
$$
and consequently $x_1^{a_1-1}\cdots x_n^{a_n-1}t^{b-g}\in\overline{R[It]}$, 
a contradiction to the choice of $x^at^b$. Thus $b\leq g$. Using the
third condition we get  
$x_1\cdots x_n\in I^g\subset I^b$, which readily implies 
$x^at^b\in R[It]$. \QED

\medskip

According to Corollary~\ref{dec11-03} condition ($\mathrm{a}_3$) is
redundant when $I$ is generated by monomials of the same degree.

\begin{Proposition} Let $Y\subset X$ and let $I_{Y}=I\cap K[Y]$. If
$I_Y$ has  
the K\"onig property for all $Y$ and $\overline{R[It]}$ is generated 
as a $K$-algebra by monomials 
of the form $x^at^b$, with $x^{a}$ square-free, then $R[It]$ is normal. 
\end{Proposition}

\demo Take $x^at^b$ a generator of $\overline{R[It]}$, with $x^a$
square-free. By  
induction we may assume $x^at^b=x_1\cdots x_nt^b$. Hence, since
$(1,\ldots,1,b)$  is in $\mathbb{R}_+{\cal A}'$, we get that 
$|C_k|\geq b$ for $k=1,\ldots,s$. 
In particular $g={\rm ht}(I)\geq b$. As $I$ has the K\"onig property,
we get  
$x_1\cdots x_n\in I^g$ and consequently $x^at^b\in R[It]$. \QED

\begin{Proposition} Let $I_i=I\cap K[X\setminus\{x_i\}]$. If $I_i$ is
normal for 
$i=1,\ldots,n$ and
\begin{equation}\label{jan6-04}
C=H_{\ell_1}\cap H_{\ell_2}\cap\cdots\cap H_{\ell_r}\cap
\mathbb{R}_{+}^{n+1}\neq(0),
\end{equation}
then $R[It]$ is normal.
\end{Proposition}

\demo Let $x^at^b=x_1^{a_1}\cdots x_n^{a_n}t^b\in\overline{R[It]}$ be
a minimal generator,  
that is $(a,b)$, cannot be written as a sum of two non-zero integral
vectors in 
$\mathbb{R}_+{\cal A}'$. It suffices to prove that $0\leq b\leq 1$
because this readily implies that $x^a$ is either a variable or a 
monomial in $F$. Assume $b\geq 2$. Since $I_i$ is normal we may
assume that $a_i\geq 1$ for all $i$. As each variable occurs in at
least one monomial of $F$, using that $C$ is contained in
$\mathbb{R}_+{\cal A}'$ together with Eq.~(\ref{jan6-04}), it follows that
there is $(v_k,1)$ such that 
$\langle (v_k,1),\ell_i\rangle=0$ for $i=1,\ldots,r$. Therefore 
$$
\langle(a-v_k,b-1),\ell_i\rangle\geq 0\ \ \ \ \ (i=1,\ldots,r).
$$
Thus $(a,b)-(v_k,1)\in \mathbb{R}_+{\cal A}'$, a contradiction to the
choice of $x^at^b$.\QED

\section{Some applications to Rees algebras and clutters}

Throughout this section we assume $\deg(x^{v_i})=d\geq 2$ for all $i$. 
By assigning $\deg(x_i)=1$ and 
$\deg(t)=-(d-1)$, the Rees algebra $R[It]$ becomes a 
standard graded $K$-algebra, i.e., it is generated by elements of
degree $1$. The $a$-invariant of $R[It]$, with 
respect to this grading, is denoted by $a(R[It])$. If $R[It]$ is a normal  
domain, then according to a formula of 
Danilov-Stanley \cite[Theorem~6.3.5]{BHer} its {\it 
canonical module\/} is the ideal 
of $R[It]$ given by
$$
\omega_{R[It]}=(\{x_1^{a_1}\cdots x_n^{a_n}
t^{a_{n+1}}\vert\, a=(a_i)
\in(\mathbb{R}_+{\cal A}')^{\rm o}\cap\mathbb{Z}^{n+1}\}),
$$
where $(\mathbb{R}_+{\cal A}')^{\rm o}$ is the topological 
interior of the Rees cone.

\begin{Theorem}\label{nov8-03} If ${\rm gr}_I(R)$ is reduced, then 
$$
a(R[It])\geq -\left[n-(d-1)(\alpha_0({\cal C})-1)\right],
$$
with equality if $I$ is unmixed. 
\end{Theorem}

\demo It is well known (see
\cite{BHer})  that the $a$-invariant can be expressed as
$$
a(R[It])=-{\rm min}\{\, i\, \vert\, 
(\omega_{R[It]})_i\neq 0\}.
$$
Set $\alpha_0=\alpha_0({\cal C})$. Using Eq.~(\ref{okayama-oct}) it
is seen that the vector  
$(1,\ldots,1,\alpha_0-1)$ is in the interior of the Rees cone. Thus
the inequality  
follows by computing the degree of $x_1\cdots x_nt^{\alpha_0-1}$.

Assume that $I$ is unmixed. Take an arbitrary monomial 
$x^at^b=x_1^{a_1}\cdots x_n^{a_n}t^b$ in the ideal $\omega_{R[It]}$, 
that is, $(a,b)\in(\mathbb{R}_+{\cal A}')^{\rm o}$. By 
Proposition~\ref{oct30-03} the vector $(a,b)$ has positive entries and 
satisfies
\begin{eqnarray}
&-b+\sum_{x_i\in C_k}a_i\geq 1 \ \ \ \ (k=1,\ldots,s).&\label{nov9-03}
\end{eqnarray}
If $\alpha_0\geq b+1$, we obtain the inequality 
\begin{equation}\label{jan9-04}
\deg(x^at^b)=a_1+\cdots+a_n-b(d-1)\geq n-(d-1)(\alpha_0-1).
\end{equation}
Now assume $\alpha_0\leq b$. Using the normality of $R[It]$ and
Eqs.~(\ref{okayama-oct})  
and (\ref{nov9-03})
it follows that the monomial
$$
m=x_1^{a_1-1}\cdots x_n^{a_n-1}t^{b-\alpha_0+1}
$$
belongs to $R[It]$. Since $x^at^b=mx_1\cdots x_nt^{\alpha_0-1}$, 
the inequality (\ref{jan9-04}) also holds in this case. Altogether we
conclude the desired  
equality. \QED

\begin{Corollary}[\cite{alexdual}]\label{nov9-1-03} If $I$ is 
unmixed with $\alpha_0({\cal C})=2$ and 
${\rm gr}_I(R)$ is reduced, then $R[It]$ is a Gorenstein ring and 
$$
a(R[It])=-(n-d+1).
$$
\end{Corollary}

\demo From the proof of Theorem~\ref{nov8-03} it follows that $x_1\cdots x_nt$ 
generates the canonical module. \QED

\medskip

Notice that if $\alpha_0({\cal C})\geq 3$, then $R[It]$ is not
Gorenstein because  
the monomials $x_1\cdots x_nt^{\alpha_0-1}$ and $x_1\cdots x_nt$ 
are distinct minimal generators of $\omega_{R[It]}$. This holds in a
more general  
setting (see Proposition~\ref{nov9-3-03} below).

\begin{Corollary} Let $J=I_c({\cal C})$ be the ideal
of vertex covers of $\cal C$. If  
$\cal C$ is a bipartite graph and $I=I({\cal C})$ is unmixed, then $R[Jt]$ 
is a Gorenstein ring and 
$$
a(R[Jt])=-(n-\alpha_0({\cal C})+1).
$$
\end{Corollary}

\demo Notice that $R[Jt]$ has the grading induced by assigning $\deg(x_i)=1$ and 
$\deg(t)=1-\alpha_0({\cal C})$. Thus the formula 
follows from Corollary~\ref{nov9-1-03} once we recall that 
${\rm gr}_J(R)$ is a reduced ring according to Corollary~\ref{nov9-2-03}. \QED

\begin{Lemma}{\rm(\cite{Bd}, \cite[p.~142]{Vas})}\label{jan10-03}
If $S$ is a regular local ring and $J$ is an ideal of $S$ generated 
by a regular sequence 
$h_1,...,h_g$, then $S[Jt]$ is determinantal: 
$$
S[Jt]\simeq S[z_1,\ldots,z_g]/I_2
\left(\hspace{-1mm}
\begin{array}{ccc}
z_1&\cdots&z_g\\ 
h_1&\cdots&h_g
\end{array}\hspace{-1mm}
\right)
$$ 
and its canonical module is $\omega_S(1,t)^{g-2}$. 
\end{Lemma}

\begin{Proposition}\label{nov9-3-03}
If $I$ has height $g\geq 2$ and $S=R[It]$ is Gorenstein, then $g=2$.
\end{Proposition}

\demo Since $I_{\mathfrak{p}}$ is a complete intersection for all associated prime 
ideals $\mathfrak{p}$ of $I$ and $S$ is Gorenstein one has 
$\omega_S\simeq\omega_R(1,t)^{g-2}$ \cite{Her}. Then
\begin{equation}\label{may5-06}
S\simeq \omega_S\simeq\omega_R(1,t)^{g-2}=R\oplus Rt
\oplus\cdots\oplus Rt^{g-2}\oplus It^{g-1}\oplus\cdots
\end{equation}
Take a minimal prime $\mathfrak{p}$ of $I$ of height $g$. Then 
$S_\mathfrak{p}=R_\mathfrak{p}[I_\mathfrak{p}t]$ is the Rees algebra 
of the ideal $I_\mathfrak{p}$, which is generated by a regular sequence. 
Thus localizing the extremes of Eq.~(\ref{may5-06}) at $\mathfrak{p}$ and 
using Lemma~\ref{jan10-03} 
we obtain 
$$
S_\mathfrak{p}=R_\mathfrak{p}[I_\mathfrak{p}t]\simeq
\omega_{R_\mathfrak{p}}(1,t)^{g-2} \simeq \omega_{S_\mathfrak{p}}. 
$$
Note that it is important to know a priori that the canonical 
module of $S_\mathfrak{p}$ is  $\omega_{R_\mathfrak{p}}(1,t)^{g-2}$.
Hence $S_\mathfrak{p}$ is  
Gorenstein. To finish the proof note that the only 
Gorenstein determinantal rings that occur in 
Lemma~\ref{jan10-03} are those with $g=2$. Here the hypothesis on the degrees 
of $x^{v_i}$ is not needed. \QED

\begin{Lemma}\label{nov19-03} If $\overline{R[It]}=R_s(I)$, 
then there is a minimal vertex cover $C_k$ 
of $\cal C$ such that $|{\rm supp}(x^{v_i})\cap C_k|=1$ for $i=1,\ldots,q$.  
\end{Lemma}

\demo We claim that $J_k=\mathfrak{p}_kR[It]$ for some $1\leq k\leq s$. 
If not, using Eq.~(\ref{may2-06}), we can pick $x^{v_k}t\in J_k$ for
$k=1,\ldots,s$. Then by  
Proposition~\ref{oct28-03} the product of these monomials is 
in the radical of $IR[It]$. Therefore 
$$
\left[(x^{v_1}t)\cdots(x^{v_s}t)\right]^p\in IR[It]
$$
for some $0\neq p\in\mathbb{N}$. Thus 
$(x^{v_1}\cdots x^{v_s})^p\in I^{sp+1}$. By degree 
considerations, using that $\deg(x^{v_i})=d$ for all $i$, one 
readily derives a contradiction. This proves the claim. Hence 
$\langle (v_i,1),\ell_k\rangle=0$ for all $i$ and 
${v_1},\ldots,{v_q}$ lie on the hyperplane 
$$
\textstyle\sum_{x_i\in C_k}x_i=1.
$$
Therefore $|{\rm supp}(x^{v_i})\cap C_k|=1$ for all $i$, as 
desired. \QED

\begin{Proposition} If $\overline{R[It]}=R_s(I)$ and $I$ is unmixed,
then
\[
H_{\ell_1}\cap H_{\ell_2}\cap\cdots\cap H_{\ell_r}\cap
\mathbb{R}_{+}^{n+1}\neq(0).
\]
\end{Proposition}

\demo Let $J=I_c({\cal C})$ be the Alexander dual of $I$. Using 
Proposition~\ref{nov3-03} one has $\overline{R[Jt]}=R_s(J)$. Thus by 
Lemma~\ref{nov19-03} there is $v_k$ such that 
$|{\rm supp}(x^{v_k})\cap C_i|=1$ for $i=1,\ldots,r$. This means that 
$(v_{k},1)$ is in the intersection of $H_{\ell_1},\ldots,H_{\ell _r}$.
\QED

\begin{Proposition}\label{nov20-03} If $\overline{R[It]}=R_s(I)$, 
then 
there are $C_1,\ldots,C_d$ mutually disjoint minimal vertex covers of $\cal C$ such 
that $\cup_{i=1}^q{\rm supp}(x^{v_i})=\cup_{i=1}^d C_i$  and 
$$
|{\rm supp}(x^{v_i})\cap C_k|=1\ \ \ \ \forall\ i,k.
$$
\end{Proposition}

\demo By induction on $d$. By Lemma~\ref{nov19-03} there is a minimal vertex 
cover $C_1$ of $\cal C$ such that $|{\rm supp}(x^{v_i})\cap C_1|=1$ for all 
$i$. Consider the ideal $I'$ obtained from $I$ by making $x_i=1$ for
all $x_i\in C_1$. Then $I'$ is an ideal generated by monomials of
degree $d-1$ and  $\overline{R[I't]}=R_s(I')$ by 
Corollary~\ref{nov19-1-03}. Thus we can apply 
induction to get the required assertion. \QED

\begin{Corollary}\label{dec11-03}
If $I$ is unmixed and $\overline{R[It]}=R_s(I)$, then both $\cal C$ and
the clutter  
$\cal D$ of minimal vertex covers of $\cal C$ have the K\"onig property. 
\end{Corollary}

\demo That $\cal D$ has the K\"onig property follows from
Proposition~\ref{nov20-03}, because $\alpha_0({\cal D})=d$ and
$C_1,\ldots,C_d$ are independent edges of $\cal D$. Now $I_c({\cal
C})$ is unmixed, is generated by 
monomials of degree $\alpha_0({\cal C})$, and according to
Proposition~\ref{nov3-03} one 
has $\overline{R[I_c({\cal C})]}=R_s(I_c({\cal C}))$. Thus, using
again Proposition~\ref{nov20-03}, we conclude that $\cal C$ has the K\"onig property. 
\QED

\medskip

Combining Corollary~\ref{dec11-03} with Proposition~\ref{nov20-1-03}
we obtain: 

\begin{Theorem}\label{may4-06} Let $I_i=I\cap K[X\setminus\{x_i\}]$. If $I$ is
unmixed and $Q(A)$ 
is integral, then ${\rm gr}_I(R)$ is reduced if and only if $I_i$ is normal 
for $i=1,\ldots,n$. 
\end{Theorem}

\bibliographystyle{plain}

\end{document}